\documentclass[11pt]{amsart}

\usepackage[dvipsnames]{xcolor}
\usepackage[pagebackref,backref=true,colorlinks=true, linkcolor=Sepia, citecolor=Sepia]{hyperref}
\usepackage{amsmath,amsthm,amssymb,fancyhdr,graphicx,bbm,cancel,mathrsfs,todonotes,xypic,pinlabel}
\usepackage{tikz}
\usepackage{extarrows,tipa}
\usetikzlibrary{matrix}
\usepackage[all]{xy}
\usepackage{mathtools}
\usepackage{bm}
\usepackage{verbatim}
\usepackage{microtype}
\usepackage{subcaption}

\theoremstyle{theorem}
\newtheorem{thm}{Theorem}[section]

\newtheorem{theoremalpha}{Theorem}

\newtheorem{lem}[thm]{Lemma}

\newtheorem{cor}[thm]{Corollary}
\theoremstyle{definition} 

\newtheorem{rmk}[thm]{Remark}
\theoremstyle{remark} 

\usepackage{enumerate,pdfpages,stmaryrd}
\usepackage[margin=1.2in, marginparwidth=1in]{geometry}
\usepackage{tikz}
\usepackage{dsfont}
\usetikzlibrary{matrix}

\renewcommand{\xto}{\xrightarrow}

\newcommand{\C}{\mathbb{C}}

\newcommand{\Z}{\mathbb{Z}}
\newcommand{\Q}{\mathbb{Q}}
\newcommand{\R}{\mathbb{R}}

\newcommand{\beq}{\begin{equation*}}
\newcommand{\eeq}{\end{equation*}}

\newcommand{\SL}{\mathrm{SL}}
\newcommand{\GL}{\mathrm{GL}}
\newcommand{\onto}{\twoheadrightarrow}

\DeclareMathOperator{\Diff}{Diff}

\DeclareMathOperator{\Hom}{Hom}

\DeclareMathOperator{\Aut}{Aut}

\DeclareMathOperator{\SO}{SO}

\DeclareMathOperator{\Out}{Out}

\newcommand{\dmo}{\DeclareMathOperator}

\newcommand{\Om}{\Omega}

\newcommand{\wtil}{\widetilde}

\newcommand{\bb}[1]{\mathbb{#1}}

\dmo{\sgn}{sign}\dmo{\Span}{span}
\dmo{\we}{\wedge}
\dmo{\ind}{ind}\dmo{\Ind}{Ind}
\dmo{\bop}{\bigoplus}\dmo{\pic}{Pic}
\dmo{\vol}{Vol}\dmo{\gal}{Gal}\dmo{\perm}{Perm}
\dmo{\tor}{Tor}\dmo{\ext}{Ext}\dmo{\Ext}{Ext}
\dmo{\aut}{aut}

\dmo{\inn}{Inn}\dmo{\var}{Var}
\dmo{\ad}{ad}\dmo{\curl}{curl}
\dmo{\hy}{\bb H}\dmo{\Sl}{SL}
\dmo{\psl}{PSL}
\dmo{\iso}{iso}
\dmo{\conf}{Conf}
\dmo{\stab}{Stab}\dmo{\Jac}{Jac }
\dmo{\diam}{diam}\dmo{\fix}{Fixed}\dmo{\Fix}{Fix}
\dmo{\injR}{injRad}\dmo{\Ad}{Ad}
\dmo{\esv}{ess-vol}
\dmo{\nil}{Nil}\dmo{\sol}{Sol}
\dmo{\Div}{div}
\dmo{\SU}{SU}
\dmo{\rk}{rk}
\dmo{\rank}{rank}
\dmo{\psp}{PSp}\dmo{\psu}{PSU}
\dmo{\PU}{PU}\dmo{\pgl}{PGL}
\dmo{\Mod}{Mod}\dmo{\range}{Range}
\dmo{\eu}{eu}\dmo{\mi}{mi}
\dmo{\Log}{Log}\dmo{\supp}{supp}
\dmo{\maps}{Maps}\dmo{\Gr}{Gr}
\dmo{\Pin}{Pin}
\dmo{\Spin}{Spin}\dmo{\Str}{Str}
\dmo{\Sq}{Sq}\dmo{\Symp}{Symp}
\dmo{\pd}{PD}\dmo{\PD}{PD}\dmo{\sig}{Sig}
\dmo{\ev}{ev}\dmo{\St}{St}
\dmo{\Pt}{Pt}\dmo{\pt}{pt}

\dmo{\Pl}{PL}
\dmo{\String}{String}\dmo{\smear}{smear}

\dmo{\dev}{dev}
\dmo{\met}{Met}\dmo{\contact}{Contact}
\dmo{\teich}{Teich}\dmo{\Teich}{Teich}\dmo{\qi}{QI}
\dmo{\der}{Der}
\dmo{\cl}{Cliff}\dmo{\Cl}{Cl}
\dmo{\Pf}{Pf}
\dmo{\ch}{ch}\dmo{\diag}{diag}
\dmo{\grad}{grad}\dmo{\Char}{char}
\dmo{\spec}{Spec}\dmo{\Arg}{Arg}
\dmo{\gl}{GL}
\dmo{\sym}{Sym}\dmo{\Sym}{Sym}
\dmo{\com}{Comm}
\dmo{\Lk}{Lk}
\dmo{\CAT}{CAT}
\dmo{\Rep}{Rep}
\dmo{\Res}{Res}
\dmo{\Conf}{Conf}
\dmo{\PConf}{PConf}
\dmo{\Push}{Push}
\dmo{\Cont}{Cont}
\dmo{\sm}{\setminus}
\dmo{\vn}{\varnothing}
\dmo{\disk}{\mathbb D}
\dmo{\Trd}{Trd}\dmo{\Mat}{Mat}
\dmo{\Riem}{Riem}
\dmo{\Diffn}{\Diff_0}\dmo{\diff}{diff}
\dmo{\homeo}{Homeo}

\dmo{\Ham}{Ham}\dmo{\Met}{Met}
\dmo{\Ein}{Ein}\dmo{\CP}{\co P}
\dmo{\Per}{Per}\dmo{\Ric}{Ric}
\dmo{\Nrd}{Nrd}
\dmo{\Comp}{Comp}\dmo{\PSC}{PSC}
\dmo{\Cent}{Cent}\dmo{\Orb}{Orb}
\dmo{\aind}{a-ind}\dmo{\tind}{t-ind}
\dmo{\constant}{constant}
\dmo{\Td}{Td}
\dmo{\LMod}{LMod}
\dmo{\SMod}{SMod}
\dmo{\SDiff}{SDiff}
\dmo{\Br}{Br}
\dmo{\csch}{csch}
\dmo{\triv}{triv}
\dmo{\genus}{genus}
\dmo{\Homeq}{HomEq}
\dmo{\PP}{\mathbb{P}}
\dmo{\U}{U}
\dmo{\Gal}{Gal}
\dmo{\BDiff}{\wtil{\Diff}}
\dmo{\BAut}{\wtil{\Aut}}
\dmo{\Iso}{Iso}
\dmo{\Cone}{Cone}
\dmo{\codim}{codim}
\dmo{\II}{II}
\dmo{\I}{I}
\dmo{\InjRad}{InjRad}
\dmo{\Inn}{Inn}
\dmo{\sys}{sys}
\dmo{\Comm}{Comm}
\dmo{\PO}{PO}
\dmo{\vertex}{Vert}
\dmo{\POm}{P\Om}
\dmo{\ab}{ab}
\dmo{\PSO}{PSO}
\dmo{\CRS}{CRS}
\dmo{\Diffext}{Diffext}
\dmo{\Diffextad}{Diffextad}
\dmo{\Diffstand}{Diffstand}

\begin{document}
\title[Rotation index, MMN pairing, and group actions]{Rotation index, Milnor--Munkres--Novikov pairing, and group actions on manifolds}
\author{Mauricio Bustamante and Bena Tshishiku}

\maketitle

\begin{abstract}
We introduce an invariant of a pair of commuting invertible matrices that we call the rotation index. We apply this invariant, together with the Milnor--Munkres--Novikov pairing, to the study of some questions about group actions of $\Z^2$, specifically the Nielsen realization problem, higher-rank Anosov actions, and extending actions from the sphere $S^{d-1}$ to the disk $D^d$. 
\end{abstract}
\section{Introduction}
\subsection{Rotation index}

We begin by introducing the invariant whose study and application is the focus of this paper. Given a pair $A,B\in\SL_d(\R)$, we define the \emph{rotation index}  as 
\begin{equation}\label{eqn:rotation-index}
\rho(A,B)=\sum_{\lambda,\mu}\dim \big[E_\lambda(A)\cap E_\mu(B)
\big]\pmod 2,
\end{equation}
where the sum ranges over pairs $\lambda,\mu$ of negative real eigenvalues of $A,B$, respectively, and $E_\lambda(A)=\ker[(A-\lambda I)^m]$ denotes the generalized $\lambda$-eigenspace for $A$.

While $\rho(A,B)$ is defined for any pair $A,B$, it is most interesting as an invariant of commuting pairs. The space of commuting pairs is defined as the space $\Hom\big(\Z^2,\SL_d(\R)\big)$ of group homomorphisms from $\Z^2$ to $\SL_d(\R)$, topologized as a subspace of $\SL_d(\R)\times\SL_d(\R)$.

\begin{theoremalpha}\label{thm:RI-continuous}
For all $d\geq 2$, the function $\rho:\SL_d(\R)\times\SL_d(\R)\to\Z/2\Z$ is continuous when restricted to the subspace of commuting pairs.
\end{theoremalpha}

By Theorem \ref{thm:RI-continuous}, $\rho$ is an invariant of path components of $\Hom\big(\Z^2,\SL_d(\R)\big)$. It is known by special cases of work of Pettet--Souto \cite{PS} and Higuera Rojo \cite{rojo} that $\Hom\big(\Z^2,\SL_d(\R)\big)$ has two components when $d\geq 3$. We give an independent proof of this using the rotation index.
\begin{theoremalpha}\label{thm:path-components}
For all $d\geq 3$, the rotation index defines a bijection 
\beq
\pi_0\big[\Hom\big(\Z^2,\SL_d(\R)\big)\big]\to\Z/2\Z.
\eeq
\end{theoremalpha}
The rotation index is related to the second Stiefel--Whitney class $w_2(\phi)\in H^2(\Z^2;\Z/2\Z)$ of a representation $\phi:\Z^2\to\SL_d(\R)$, $d\geq 3$, which is the obstruction to lifting $\phi$ to the universal cover $\wtil{\SL_d(\R)}\to\SL_d(\R)$. Identifying $H^2(\Z^2;\Z/2\Z)\cong H^2(T^2;\Z/2\Z)$, we may evaluate on a fundamental class $[T^2]\in H_2(T^2;\Z/2\Z)$ to obtain $\langle w_2(\phi),[T^2]\rangle\in\Z/2\Z$. 
Since path components of $\Hom\big(\Z^2,\SL_d(\R)\big)$ correspond to values of the Stiefel--Whitney class, and hence to values of $\rho$, Theorem \ref{thm:path-components} implies the following corollary.
\begin{cor}
\label{cor:steifel-whitney}Fix $d\ge3$. For any representation $\phi:\Z^2\to\SL_d(\R)$, we have 
\beq
\langle w_2(\phi),[T^2]\rangle=\rho(\phi).
\eeq
\end{cor}
Here we write $\rho(\phi)$ to mean $\rho(A,B)$, where $(A,B)=(\phi(1,0),\phi(0,1))$ is a commuting pair.

In the remainder of the introduction we illustrate how to apply the rotation index $\rho$ to study group actions. For each of our applications, we will also use the Milnor--Munkres--Novikov pairing (see \cite[Section 1.3.2]{BKKT} and references therein)
\[\pi_1^s\times\Theta_d\to\Theta_{d+1}\]
where $\Theta_d$ denotes the Kervaire--Milnor group of homotopy $d$-spheres \cite{kervaire-milnor} and $\pi_1^s\cong\Z/2\Z$ is the first stable homotopy group of spheres. Below we write $\eta\in\pi_1^s\cong\Z/2\Z$ for the generator, and denote the pairing with $\eta$ as a map
\[
\begin{array}{rccl}\eta:&\Theta_d&\to&\Theta_{d+1}\\[2mm]
&\Sigma&\mapsto&\eta\cdot\Sigma
\end{array}
\]

An important role will be played by the function 
\[\delta(\Sigma)=\begin{cases}0&\text{ if }2\mid\eta\cdot\Sigma\\1&\text{ if }2\nmid\eta\cdot\Sigma.\end{cases}
\]
for $\Sigma\in\Theta_d$. The invariant $\delta$ takes each of the values $0$ and $1$ infinitely often (since $\Theta_d$ is finite for each $d$, necessarily this requires taking $d$ to infinity). See \cite[Table 1 and Rmk.\ 1.10]{BKKT}.
\subsection{Homotopy Anosov actions on exotic tori} \label{sec:intro-anosov}

Recall that an exotic $d$-torus is a closed smooth manifold that is homeomorphic but not diffeomorphic to the standard torus $T^d=\R^d/\Z^d$. For example, $T^d\#\Sigma^d$ and $(T^{d-1}\#\Sigma^{d-1})\times S^1$ are exotic tori when $\Sigma$ is an exotic sphere.

It is known that some exotic tori admit Anosov diffeomorphisms \cite{FJ-anosov, farrell-gogolev} (i.e.\ Anosov actions of $\Z$), but there is no known example of an Anosov action of $\Z^2$ on an exotic torus $T^d\#\Sigma$. An Anosov action of $\Z^2$ on $T^d\#\Sigma$ determines (via the induced action on $\pi_1$) an Anosov action on the standard torus $T^d$, which is called its linearization. Thus in the study of Anosov actions on exotic tori, a basic problem is to determine when an Anosov action on $T^d$ is the linearization of an Anosov action on $T^d\#\Sigma$. For example, by work of Rodriguez-Hertz--Wang, an Anosov action on $T^d$ that is \emph{without rank-1 factor} (c.f.\ \cite[Defn.\ 2.8]{RHW}) is not the linearization of an Anosov $\Z^2$ actions on any exotic torus \cite[Cor.\ 1.2]{RHW}. Using our Theorem \ref{thm:RI-continuous}, we give examples with a stronger conclusion, where ``Anosov action" is replaced by ``smooth action".

\begin{theoremalpha}\label{thm:anosov}
For each exotic sphere $\Sigma\in\Theta_d$ such that $\delta(\Sigma)=1$, there exists an Anosov action $\Z^2\curvearrowright T^d$ without rank-$1$ factor that is not the linearization of any smooth $\Z^2$ action on $T^d\#\Sigma$. In fact, there are infinitely many conjugacy classes of such actions.
\end{theoremalpha}

For example, Theorem \ref{thm:anosov} applies for $d=8$ to the unique exotic torus of the form $T^8\#\Sigma$ (it is unique because $\Theta_8=\Z/2\Z$). 
We prove Theorem \ref{thm:anosov} in Section \ref{sec:anosov}.

\subsection{Nielsen realization and Borel's theme}

En route to Theorem \ref{thm:anosov} we prove the following stronger result.
\begin{theoremalpha}\label{thm:nielsen}
Fix $d\geq 7$, and fix an exotic sphere $\Sigma\in\Theta_d$. 
Let $G=\langle A,B\rangle<\SL_d(\Z)$ be a subgroup isomorphic to $\Z^2$. If $\delta(\Sigma)=1$, then the natural homomorphism
\[\ell:\Diff^+(T^d\#\Sigma)\to\Out^+\big(\pi_1(T^d\#\Sigma)\big)\cong \SL_d(\Z)\]
splits over $G$ if and only if $\rho(A,B)=0$.
\end{theoremalpha}

Theorem \ref{thm:nielsen} is a result about Nielsen realization, which asks when a group of symmetries of $\pi_1(M)$ can be realized by a group of symmetries of $M$. For aspherical manifolds, this question fits within ``Borel's theme" (c.f.\ Weinberger \cite{weinberger}) that asks to what extent geometric or topological properties of an aspherical manifold $M$ can be determined by the fundamental group $\pi_1(M)$. Indeed, in the setting of Theorem \ref{thm:nielsen} one can determine if $\Z^2<\SL_d(\Z)$ is realized by a group of diffeomorphisms of $T^d\#\Sigma$ just using the rotation index. 
\begin{rmk}
In Theorem \ref{thm:nielsen}, the assumption $\delta(\Sigma)=1$ implies that the surjection $\ell:\Diff^+(T^d\#\Sigma)\onto\SL_d(\Z)$ does not split by work of Krannich--Kupers and the authors \cite[Thm.\ A]{BKKT}.  
Whether $\ell$ splits over a given subgroup $\Z^2<\SL_d(\Z)$ is a more subtle question, which is answered by Theorem \ref{thm:nielsen}. 
\end{rmk}
\subsection{Bordism of group actions}

We give one more (somewhat different) application of the rotation index and Milnor--Munkres--Novikov pairing. 

\begin{theoremalpha}\label{thm:bordism}
Fix $d$ such that there exists $\Sigma\in\Theta_d$ with $\delta(\Sigma)=1$. Then there exists a smooth action of $\Z^2$ on the sphere $S^{d-1}$ by diffeomorphisms isotopic to the identity such that the action does not extend to an action of $\Z^2$ on $D^d$. 
\end{theoremalpha}

Theorem \ref{thm:bordism} fits within the study of bordism of group actions; see \cite{MN1} and references therein. The assumption in Theorem \ref{thm:bordism} that the diffeomorphisms are isotopic to the identity means that each individual diffeomorphism extends from $S^{d-1}$ to $D^d$, and yet for the examples in Theorem \ref{thm:bordism}, there is no way to extend the group action. 

The construction of the examples $\Z^2\curvearrowright S^{d-1}$ in Theorem \ref{thm:bordism} is explicit. It would be interesting to determine whether these actions can extend to any $d$-manifold  with boundary $S^{d-1}$. 

\subsection*{Proof techniques}
Our proofs of Theorems \ref{thm:RI-continuous} and \ref{thm:path-components} rely on perturbation theory of linear operators. In Theorem \ref{thm:nielsen}, the obstruction to lifting is detected on the mapping class group; here we use computations from \cite{BKKT} and an Euler class argument; computing this Euler class relates to the rotation index. Theorem \ref{thm:anosov} is deduced from Theorem \ref{thm:nielsen}, together with a simple number-theoretic construction to produce commuting hyperbolic matrices in $\SL_d(\Z)$ with certain properties. Finally, Theorem \ref{thm:bordism} is proved by identifying the Milnor-Munkres-Novikov pairing with the obstruction to split a certain central group extension.
\subsection*{Section Outline}
We establish the continuity of the rotation index on commuting matrices and prove Theorems \ref{thm:RI-continuous} and \ref{thm:path-components} in Section \ref{sec:rotation-index}.  Sections \ref{sec:proof} and \ref{sec:anosov} contain the proofs of Theorem \ref{thm:nielsen}
and \ref{thm:anosov} respectively.

\subsection*{Acknowledgements} We thank Andrey Gogolev for asking us a question that motivated our main result and thank Sebasti\'an Hurtado for useful comments. MB is supported by  ANID Fondecyt Regular grant 1250727. BT is supported by NSF grant DMS-2104346.

\section{The rotation index } \label{sec:rotation-index}

In this section we prove some basic properties about the rotation index, including Theorem \ref{thm:RI-continuous}. We also use this point-of-view to give a new proof that $\Hom(\Z^2,\SL_d(\R))$ has two components for $d\ge3$. One of the tools we use to prove Theorem \ref{thm:RI-continuous} is perturbation theory for linear operators, which we recall in \S\ref{sec:RI-background} along with some other auxiliary results.

\subsection{Perturbation theory and eigenspaces}\label{sec:RI-background}

Generalized eigenspaces do not generally vary continuously on $\SL_d(\R)$ (e.g.\ because the dimensions may change). Perturbation theory for linear operators gives some information about how the eigenspaces vary and gives a weaker form of continuity, which is helpful for our purposes. We recall this below. Our reference is \cite[Ch.\ 2, \S1]{kato}.

Fix $A\in\SL_d(\R)$. Let $\epsilon>0$ be a small real number, and denote by $T^1_I\SL_d(\R)$ the subspace of unit vectors in the tangent space of $\SL_d(\R)$ at the identity $I$, which we identify with the space of traceless real $d\times d$ matrices. We consider perturbations of $A$ of the form
\beq
A(t,V)=A\cdot\exp(tV),
\eeq
where $t\in [0,\epsilon)$, $V\in T^1_I\SL_d(\R)$, and $\exp(-)$ denotes the matrix exponential. Observe that $A(t,V)$ is a continuous function, satisfies $A(0,V)=A$, and lies in $\SL_d(\R)$ for all $(t,V)\in [0,\epsilon)\times T^1_I\SL_d(\R)$. 

Recall that the \textit{resolvent}
\beq
R(\zeta,t,V)=(A(t,V)-\zeta)^{-1}
\eeq
of $A(t,V)$ is defined for all complex numbers $\zeta$ which are not eigenvalues of $A(t,V)$.
The \emph{resolvent set} of $A$ is the set $\C\setminus\spec(A)$, where $\spec(A)$ is the set of eigenvalues of $A$. 

\begin{lem}\label{lem:resolvent}
Let $A\in\SL_d(\R)$ and
\beq
\mathcal U\subset \C\times[0,\varepsilon)\times T^1_I\SL_d(\R)
\eeq
be an open set with the property that for every $(\zeta,t,V)\in\mathcal U$ the number $\zeta$ is not an eigenvalue of $A(t,V)$. Then the resolvent
\beq
R(\zeta,t,V)=(A(t,V)-\zeta)^{-1}
\eeq
is continuous on $\mathcal U$. Moreover $R(\zeta,t,V)$ tends to $(A-\zeta)^{-1}$ as $t\to 0$ uniformly for $\zeta$ ranging over compact subsets of the resolvent set of $A$.
\end{lem}
\begin{proof}
Continuity on $\mathcal U$ follows from the fact that matrix inversion in $\GL_d(\C)$ is continuous.

For the second claim, fix a compact subset $K$ of the resolvent set of $A$. For $t$ sufficiently small, $R(\zeta,t,V)$ is defined for all $\zeta\in K$ and $V\in T^1_I\SL_d(\R)$. Uniform continuity of $R(\zeta,t,V)$ then follows from compactness.
\end{proof}

Now let $\lambda\in\C$ be an eigenvalue of $A$ of algebraic multiplicity $m$. Let $\Gamma$ be a simple closed curve in $\C$ enclosing $\lambda$ but no other eigenvalues of $A$. Set
\beq
\delta=\min\{||(A-\zeta)^{-1}||^{-1}\ |\ \zeta\in\Gamma\}.
\eeq 
Note that $\delta>0$ and $R(\zeta,t,V)$ exists for all $\zeta\in\Gamma$ if $||A(t,V)-A||<\delta$. Thus the following operator on $\C^d$ can be defined:
\beq
P(t,V):=-\frac{1}{2\pi i}\displaystyle\int_{\Gamma}R(\zeta,t,V)\, d\zeta.
\eeq
Lemma \ref{lem:resolvent} yields the following corollary.
\begin{cor}\label{cor:projection}
For $(t,V)$ as above, the operator $P(t,V)$ is continuous.
\end{cor}
The operator $P(t,V)$ is the projection to the eigenspaces of $A(t,V)$ with eigenvalues lying inside the curve $\Gamma$. In particular, the image of $P(t,V)$ is the sum $\bigoplus E_{\lambda'}(A(t,V))$ of the generalized eigenspaces $E_{\lambda'}(A(t,V))$ for the eigenvalues $\lambda'$ of $A(t,V)$ lying in $\Gamma$ \cite[I. Problem 5.9]{kato}.

We write $\widehat{\lambda}=\widehat{\lambda}(t,V)$ for the eigenvalues of $A(t,V)$ lying inside of $\Gamma$, and we refer to these eigenvalues as the $\lambda$-group. The sum of generalized eigenspaces ranging over all $\lambda'$ in the $\lambda$-group will be abbreviated as 
\beq
E_{\widehat{\lambda}}(A(t,V)):=\bigoplus E_{\lambda'}(A(t,V)).
\eeq

\begin{lem}\label{lem:Dim-intersection}
Let $(A,B)\in\SL_d(\R)\times\SL_d(\R)$ be a commuting pair. There exists an open neighborhood $\mathcal U\subset \SL_d(\R)\times\SL_d(\R)$ of $(A,B)$ such that for every pair $(A',B')\in\mathcal U$ with $A'B'=B'A'$ and for every $\lambda\in\spec(A)$ and $\mu\in\spec(B)$ we have
\beq
\dim\big(E_\lambda(A)\cap E_\mu(B)\big)
=
\dim\big(E_{\widehat\lambda}(A')\cap E_{\widehat\mu}(B')\big),
\eeq
where $\widehat\lambda$ (resp.\ $\widehat\mu$) denotes the $\lambda$-group for $A'$ (resp.\ the $\mu$-group for $B'$) determined by a fixed isolating contour around $\lambda$ (resp.\ $\mu$).
\end{lem}

\begin{proof}
Fix $\lambda\in\spec(A)$ and $\mu\in\spec(B)$. Choose simple closed curves $\Gamma_\lambda$ and $\Gamma_\mu$ in $\C$ such that $\Gamma_\lambda$ encloses $\lambda$ and no other eigenvalue of $A$, and $\Gamma_\mu$ encloses $\mu$ and no other eigenvalue of $B$. Define the spectral projections
\beq
P:=-\frac{1}{2\pi i}\int_{\Gamma_\lambda} (A-\zeta)^{-1}\,d\zeta,
\qquad
Q:=-\frac{1}{2\pi i}\int_{\Gamma_\mu} (B-\eta)^{-1}\,d\eta.
\eeq
As mentioned above, $P$ is the projection onto the generalized eigenspace $E_\lambda(A)$ and $Q$ is the projection onto $E_\mu(B)$. In particular $P$ and $Q$ are idempotent linear operators of finite rank.

Since $A$ and $B$ commute, so do $P$ and $Q$. Therefore $PQ$ is again a projection and its image satisfies
\beq
\operatorname{im}(PQ) = \operatorname{im}P \cap \operatorname{im}Q = E_\lambda(A)\cap E_\mu(B).
\eeq
Thus
\beq
\dim\big(E_\lambda(A)\cap E_\mu(B)\big)=\operatorname{rank}(PQ).
\eeq

Set
\beq
\delta_A := \min_{\zeta\in\Gamma_\lambda}||(A-\zeta)^{-1}||^{-1}>0,
\qquad
\delta_B := \min_{\eta\in\Gamma_\mu}||(B-\eta)^{-1}||^{-1}>0.
\eeq

By Corollary \ref{cor:projection}, whenever $||A'-A||<\delta_A$ for $A'$ of the form $A'=A\cdot\exp(tV)$, the resolvent $(A'-\zeta)^{-1}$ exists for every $\zeta\in\Gamma_\lambda$ and depends continuously on $A'$, and similarly whenever $||B'-B||<\delta_B$ the resolvent $(B'-\eta)^{-1}$ exists for every $\eta\in\Gamma_\mu$ and depends continuously on $B'=B\cdot\exp(tW)$. Hence any subset $\mathcal V\subset\SL_d(\R)\times\SL_d(\R)$ in the image of the product map
\beq
[0,\delta_A)\times T^1_I\SL_d(\R)\times [0,\delta_B)\times T^1_I\SL_d(\R)\xto{A\cdot\exp\times B\cdot\exp} \SL_d(\R)\times\SL_d(\R)
\eeq
is an open neighborhood of $(A,B)$ on which the contour integrals defining spectral projections with respect to $\Gamma_\lambda$ and $\Gamma_\mu$ are well-defined.

For $(A',B')\in\mathcal V$ define
\beq
P':=-\frac{1}{2\pi i}\int_{\Gamma_\lambda}(A'-\zeta)^{-1}\,d\zeta,\qquad
Q':=-\frac{1}{2\pi i}\int_{\Gamma_\mu}(B'-\eta)^{-1}\,d\eta.
\eeq
These are continuous functions of $(A',B')\in\mathcal V$. For such $(A',B')$ the images satisfy 
\beq
\operatorname{im}(P') = E_{\widehat\lambda}(A')\qquad \text{and}\qquad \operatorname{im}(Q') = E_{\widehat\mu}(B'),
\eeq 
where $\widehat\lambda,\widehat\mu$ denote the $\lambda$- and $\mu$- groups of eigenvalues of $A'$ and $B'$ inside the chosen contours.

If, in addition 
\beq
(A',B')\in\mathcal U':=\mathcal V\cap\Hom(\Z^2,\SL_d(\R)\times\SL_d(\R))
\eeq
(i.e.\ the pair is commuting), then $P'Q'=Q'P'$, and the same argument as before yields
\beq\label{eq:rank-P'Q'}
\dim\big(E_{\widehat\lambda}(A')\cap E_{\widehat\mu}(B')\big)=\operatorname{rank}(P'Q').
\eeq
In finite dimensions the rank of a projection equals its trace, and the trace is a continuous function of matrix entries; hence the map $\mathcal U'\to\R$ given by
\beq
(A',B')\mapsto\operatorname{tr}(P'Q')
\eeq
is continuous and takes values in $\Z$, so it is constant on a possibly smaller neighborhood $\mathcal U\subset\mathcal U'$.

Therefore, for every commuting pair $(A',B')\in\mathcal U$ we have
\beq
\dim\big(E_\lambda(A)\cap E_\mu(B)\big)
= \operatorname{rank}(PQ)
= \operatorname{rank}(P'Q')
= \dim\big(E_{\widehat\lambda}(A')\cap E_{\widehat\mu}(B')\big).
\eeq
This proves the lemma.
\end{proof}
The following fact about commuting matrices will be useful in proving Theorem \ref{thm:RI-continuous}.
Write $\psi:\GL_n(\C)\to\GL_{2n}(\R)$ for the homomorphism induced by the ring homomorphism $\C\to \mathrm{M}_2(\R)$ given by $x+iy\mapsto\left(\begin{smallmatrix}
x & y\\
-y & x
\end{smallmatrix}\right)$.

\begin{lem}\label{lem:Commute-realification}
Fix $\lambda\in\C\setminus\R$. If $B\in\SL_{2n}(\R)$ commutes with $\psi(\lambda I_n)$, then $B$ is also in the image of $\psi$.
\end{lem}

\begin{proof}
Writing $\lambda=x+iy$, we can write $\psi(\lambda I_n)=xI_{2n}+yJ_{2n}$, where $J_{2n}$ is the $(2n)\times(2n)$ block diagonal matrix with blocks $\left(\begin{smallmatrix}
0& 1\\
-1& 0
\end{smallmatrix}\right)$. Note that $y\neq0$ by our assumption that $\lambda\in\C\setminus\R$. Using this expression for $\psi(\lambda I_n)$, we observe that if $B$ commutes with $\psi(\lambda I_n)$, then $B$ also commutes with $J_{2n}$. Since $J_{2n}$ defines the given complex structure on $\R^{2n}$, the image of $\psi$ is precisely the centralizer of $J_{2n}$.
\end{proof}
\subsection*{Proof of Theorem \ref{thm:RI-continuous}}
We want to show that the rotation index is continuous on commuting pairs. Recall that $\rho(A,B)=\sum \dim\big[E_\lambda(A)\cap E_\mu(B)\big]\pmod 2$,
where the sum ranges over $\lambda\in\spec(A)$ and $\mu\in\spec(B)$ with $\lambda,\mu<0$. Fix such a pair $\lambda,\mu$.

For $(A',B')$ sufficiently close to $(A,B)$, we know from Lemma \ref{lem:Dim-intersection} that
\[\dim E_{\widehat\lambda}(A')\cap E_{\widehat\mu}(B') =\dim E_{\lambda}(A)\cap E_{\mu}(B).\]
Ultimately, we want to compare $\dim E_{\lambda}(A)\cap E_{\mu}(B)$ to the dimension of the subspace of $E_{\widehat\lambda}(A')\cap E_{\widehat\mu}(B')$ where $A',B'$ both act with negative real eigenvalues (in the sense of generalized eigenvalues).

The action of $A'$ and $B'$ on $E_{\widehat\lambda}(A')\cap E_{\widehat\mu}(B')$ is simultaneously upper triangular over $\C$ (because $A',B'$ commute). Write $d_{\lambda\mu}$ for the dimension of  $E_{\widehat\lambda}(A')\cap E_{\widehat\mu}(B')$. Fixing a basis so that $A',B'$ are upper triangular, we write $\{1,\ldots,d_{\lambda\mu}\}=N_A\sqcup C_A$ for the partition corresponding to which of the diagonal entries of $A'$ are negative or non-real complex numbers (positive eigenvalues do not occur because the eigenvalues of $A$ on this subspace are in a neighborhood of $\lambda<0$). Similarly, decompose $\{1,\ldots,d_{\lambda\mu}\}=N_B\sqcup C_B$.

With this notation, we want to show
\[d_{\lambda\mu}\equiv |N_A\cap N_B|\pmod 2.\]
Since $d_{\lambda\mu}=|N_A\cap N_B|+|N_A\cap C_B|+|C_A\cap N_B|+|C_A\cap C_B|$, it suffices to show that the sum of the last three terms is even. In fact each of these terms is even. Clearly $|N_A\cap C_B|+|C_A\cap C_B|=|C_B|$ are even since nonreal eigenvalues occur in complex conjugate pairs. Furthermore, both of $|N_A\cap C_B|$ and $|C_A\cap N_B|$ are even by Lemma \ref{lem:Commute-realification}. This completes the proof. 
\qed

We now turn to the proof of Theorem \ref{thm:path-components}. It is worth mentioning that the fact that 
\beq
\pi_0\big[\Hom\big(\Z^2,\SL_d(\R)\big)\big]\cong\Z/2\Z
\eeq for $d\ge3$ follows from work of Pettet--Souto and Higuera-Rojo \cite{PS,rojo} (Pettet--Souto show $\Hom\big(\Z^2,\SL_d(\R)\big)$ deformation retracts $\Hom\big(\Z^2,\SO(d)\big)$, and Higuera-Rojo computes components of the latter space). For our applications it will be important that the path components correspond to the values of the rotation index.
\subsection*{Proof of Theorem \ref{thm:path-components}}

The map is well-defined by Theorem \ref{thm:RI-continuous}. To see that it's surjective, observe that $\rho(I,I)=0$ and $\rho(D_1,D_2)=1$, where $D_1,D_2$ are the diagonal matrices $D_1=(-1,-1,1,\ldots,1)$ and $D_2=(-1,1,-1,1,\ldots,1)$. It remains then to show injectivity, and for this we will show that any commuting pair $(A,B)$ can be connected by a path of commuting pairs in $\SL_d(\R)$ to either $(I,I)$ or $(D_1,D_2)$. We do this in several steps.

For notational convenience we now write $(A_1,A_2)$ instead of $(A,B)$ for the commuting pair.

\subsubsection*{Path $\#1$ (to semisimple matrices).}
Any $A\in\SL_d(\R)$ can be written in the form $A=A^sA^u$ where $A^s$ is a semisimple matrix and $A^u$ is a unipotent matrix (multiplicative Jordan decomposition). Similar to \cite[Prop.\ 8.4 and its proof]{PS}, for $t\in [0,1]$ and $i=1,2$, we construct paths
\beq
A_i(t)=A_i^s\, e^{(1-t)\log A_i^u}
\eeq
where, for any square matrix $X$ and unipotent matrix $Y$, the exponential and logarithm are matrices defined by
\beq
e^X=\displaystyle\sum_{k=0}^{\infty}\frac{1}{k!}X^k\qquad\text{and}\qquad
\log Y=\displaystyle\sum_{k=1}^{\infty}\frac{(-1)^{k+1}}{k!}(Y-I)^k.
\eeq
For every $t\in [0,1]$ the matrices $A_1(t)$ and $A_2(t)$ commute. Indeed, $A_i^s$ and $A_i^u$ are polynomials in $A_i$ \cite[Ch.\ 1, Prop 4.2]{Borel}, so since $A_1,A_2$ commute, each of pair of $A_1,A_2,A_1^s,A_1^u,A_2^s,A_2^u$ commute, and similarly, $A_i^s$ and $e^{(1-t)\log A_j^u}$ commute for each pair of $i,j$ because the latter can be expressed as a power series in $A_j^u$.

At time $t=1$ of the paths $A_1(t), A_2(t)$, the matrices obtained are the semisimple parts $A_1^s,A_2^s$ of $A_1,A_2$, respectively.
\subsubsection*{Path $\#2$ (to block diagonal matrices).}
We construct commuting paths from $A_1^s, A_2^s$ to block diagonal matrices.

The matrix $A_1^s$ is semisimple, so by the rational canonical form, there exists $P\in\SL_d(\R)$ such that
\begin{equation}\label{eq:blocks}
PA_1^sP^{-1}=
\begin{psmallmatrix}
\lambda_1I_{m_1} &  &  &  \\
& \ddots  &  &  \\
 &   & \lambda_rI_{m_r} &   \\
 & &  & M_1\otimes I_{n_1} \\
   &   &  & & \ddots \\
     &   &  & && M_{s}\otimes I_{n_s}
\end{psmallmatrix}
\end{equation}
where $	\lambda_i\neq\lambda_j$ if $i\neq j$ and $M_i\in\GL_2(\R)$ is the ``realification'' of a nonreal matrix in $\GL_1(\C)$, i.e.\ each $M_i$ has the form $\left(\begin{smallmatrix}
a_i & b_i\\
-b_i & a_i
\end{smallmatrix}\right)$,
where $a_i,b_i\in\R$ and $b_i\neq 0$.

Since $A_2^s$ commutes with $A_1^s$, it preserves the eigenspaces corresponding to each real eigenvalue $\lambda_i$ and also preserves the $2n_i$-dimensional subspace on which $A_1$ acts by the matrix $M_i\otimes I_{n_i}$ (the latter subspace is the kernel of an irreducible quadratic factor of the characteristic polynomial of $A_1^s$). Consequently, we can write
\begin{equation}\label{eq:blocks2}
PA_2^sP^{-1}=
\begin{psmallmatrix}
F_1 &  &  &  \\
& \ddots  &  &  \\
 &   & F_r &   \\
 & &  & H_1 \\
  & &  & & \ddots \\
   & &  & & &H_s \\
\end{psmallmatrix}
\end{equation}
where $F_i\in\GL_{m_i}(\R)$ and $H_i\in\GL_{2n_i}(\R)$ commute with the corresponding blocks of $PA_1^sP^{-1}$.

Fix a path $P(t)$ in $\SL_d(\R)$, $t\in [0,1]$, between the identity and $P$, and consider the path
\beq
P(t)\, A_i^s\, P(t)^{-1}.
\eeq
For notational convenience, we write $B_i=PA_i^sP^{-1}$, which is the time $t=1$ endpoint of our second path.
\subsubsection*{Path $\#3$ (to diagonal matrices).}
We construct commuting paths of block diagonal matrices from $PA_1^sP^{-1},PA_2^sP^{-1}$ to matrices that are the identity on the $2n_i\times 2n_i$ blocks of both $PA_1^sP^{-1}$ and $PA_2^sP^{-1}$, c.f.\ Equations (\ref{eq:blocks}),(\ref{eq:blocks2}).

By Lemma \ref{lem:Commute-realification}, the matrices $H_i\in\GL_{2n_i}(\R)$ in Equation (\ref{eq:blocks2}) lie in the image of the realification homomorphism $\varphi:\GL_{n_i}(\C)\to\GL_{2n_i}(\R)$. 
Let $K_i\in\GL_{n_i}(\C)$ be such that $\varphi(K_i)=H_i$.
Consider a path $\alpha_1(t)$ in the subspace of scalar  matrices in $\GL_{n_i}(\C)$ that takes $\sigma_i I_{n_i}$ to the identity matrix. Any path $\alpha_2(t)$ from $K_i$ to the identity (which exists because $\GL_{n_i}(\C)$ is path-connected) will commute with the path of scalar matrices $\alpha_1(t)$. Therefore the paths $\varphi\circ\alpha_i(t)$ in $\GL_{2n_i}(\R)$, $i=1,2$, commute and take $M_i\otimes I_{n_i}$ and $H_i$ to the identity in $\GL_{2n_i}(\R)$. Combining these paths for all $i=1,2,\ldots, s$ gives the desired (commuting) paths from $PA_1^sP^{-1}$ and $PA_2^sP^{-1}$ to matrices of the form
\beq
\begin{pmatrix}
\lambda_1I_{m_1} &  &  &  \\
& \ddots  &  &  \\
 &   & \lambda_rI_{m_r} &   \\
 & &  & I_{2n} \\
\end{pmatrix}
\qquad\text{and}\qquad
\begin{pmatrix}
F_1 &  &  &  \\
& \ddots  &  &  \\
 &   & F_r &   \\
 & &  & I_{2n} \\
 \end{pmatrix}
\eeq
respectively, where $n=n_1+n_2+\cdots+n_s$.
\subsubsection*{Path $\#4$ (to diagonal matrices in $\SO(d)$).}
The sign of each $\lambda_i$ determines whether the block $\lambda_iI_{m_i}$ can be connected to either $I_{m_i}$ or $-I_{m_i}$ by a path of scalar matrices in $\GL_{m_i}(\R)$. This path will commute with any path taking $F_i$ to either $I_{m_i}$ or $R_{m_i}=\mathrm{diag}(-1,1,1,\ldots,1)$, depending on the path component of $\GL_{m_i}(\R)$ containing $F_i$. So the time $t=1$ endpoint of these paths are block diagonal matrices $B_1,B_2$ of the form
\[\mathrm{diag}(\pm I_{m_1},\ldots,\pm I_{m_r},I_{2n})\] and
\[\mathrm{diag}(C_1,\ldots, C_r,I_{2n}),\]
respectively, where $C_i\in\{I_{n_i},R_{n_i}\}$.

\subsubsection*{Path $\#5$ (to diagonal matrices in $\SO(3)$).}
Next we construct paths $\gamma_i(t)$, $0\le t\le1$, within the space of commuting matrices between $B_1$ and $B_2$ and either $(I,I)$ or the pair $(D_1,D_2)$, where
\begin{equation}\label{eqn:diagonal-pair}D_1=\mathrm{diag}(-1,-1,1,1,\ldots,1)\end{equation}and
\[D_2=\mathrm{diag}(-1,1,-1,1,\ldots,1).\]
For each of the matrices obtained in the previous step, we can pair $2\times 2$ blocks of determinant $1$  which preserve the same $2$-dimensional subspace of $\R^d$. These blocks can be joined by commuting paths in $\SO(2)$ to the identity. Along these paths, the dimension where both matrices act by multiplication by a negative scalar changes by an even number, so they can be assembled to yield paths that at time $t=1$ is equal to $(I,I)$ when $\rho(A_1,A_2)=0$ and $(D_1,D_2)$ when $\rho(A_1,A_2)=1$ (possibly after conjugating by an (even) permutation matrix, which does not change the path component).

Concatenating the preceding paths, we see that any commuting pair $(A_1,A_2)$ is connected by a path of commuting pairs in $\SL_d(\R)$ to either $(I,I)$ or $(D_1,D_2)$. This completes the proof.
\qed

We end this section with the following lemma, which is helpful in deducing Theorem \ref{thm:nielsen}.

\begin{lem}\label{lem:Subgroup}
For a commuting pair $A,B\in\SL_d(\R)$, the rotation index $\rho(A,B)$ depends only on the subgroup generated by $A,B$.
\end{lem}

\begin{proof}
We want to show that if $\langle A',B'\rangle=\langle A,B\rangle$, then $\rho(A',B')=\rho(A,B)$.

First assume that the group $\langle A,B\rangle$ is isomorphic to $\Z^2$. Any two bases for $\Z^2$ differ by an element of $\GL_2(\Z)$. The group $\GL_2(\Z)$ is generated by $\left(\begin{smallmatrix}
1 & 1\\
0 & 1
\end{smallmatrix}\right)$, $\left(\begin{smallmatrix}
1 & 0\\
1 & 1
\end{smallmatrix}\right)$, and $\left(\begin{smallmatrix}
0 & 1\\
1 & 0
\end{smallmatrix}\right)$, so it suffices to show 
\beq
\rho(A,AB)=\rho(A,B)\qquad \rho(AB,B)=\rho(A,B)\qquad \rho(B,A)=\rho(A,B).
\eeq
The last equality holds because $\rho$ is symmetric by definition. The proof of the first two are similar to each other. For definiteness, we show $\rho(A,AB)=\rho(A,B)$.

To show $\rho(A,AB)=\rho(A,B)$, first recall that since $A,B$ commute, they are simultaneously upper triangular over $\C$. Assuming $A,B$ are upper-triangular, write $\{1,\ldots,d\}=P_A\sqcup N_A\sqcup C_A$ for the partition corresponding to which of the diagonal entries of $A$ are positive, negative, or non-real complex numbers. We do the same for $B$ and for $AB$.

We want to show $|N_A\cap N_B|\equiv |N_A\cap N_{AB}|\pmod 2$. Observe that
\[N_A\cap N_{AB}=N_A\cap \big[(P_A\cap N_B)\cup (N_A\cap P_B)\big]=N_A\cap P_B,\] so we want to show that $|N_A\cap N_B|+|N_A\cap P_B|\equiv0\pmod 2$. Next observe that
\[
|N_A\cap N_B|+|N_A\cap P_B|+|N_A\cap C_B|=|N_A| \equiv 0\pmod 2.
\]
The first equality holds because $P_B\cup N_B\cup C_B=\{1,\ldots,d\}$. The second equality holds because $\det(A)=1$, so the number of its real negative eigenvalues is even. It remains to explain why $|N_A\cap C_B|$ is even. This is because $A$ and $B$ are real and commute, so the action of $A$ on the generalized $(\lambda,\bar\lambda)$-eigenspace for a complex conjugate pair of eigenvalues of $B$ is by the realification of a complex matrix, so real eigenvalues appear with even multiplicity.

The other possibilities for $\langle A,B\rangle$ are $\Z\oplus\Z/\ell\Z$ or $\Z/k\Z\oplus\Z/\ell\Z$ or $\Z$. These cases are easier because the automorphism groups are similarly generated but smaller, so the proof in the $\Z^2$ case applies the same.
\end{proof}
\section{Nielsen realization (Proof of Theorem \ref{thm:nielsen})} \label{sec:proof}

In this section we fix $\Sigma\in\Theta_d$ with $\delta(\Sigma)=1$ and we study the lifting problem
\[\ell:\Diff^+(T^d\#\Sigma)\to\Out^+\big(\pi_1(T^d\#\Sigma)\big)\cong \SL_d(\Z)\]
for subgroups of $\SL_d(\Z)$ of the form $G=\langle A,B\rangle\cong\Z^2$. To prove Theorem \ref{thm:nielsen}, we consider the cases $\rho(A,B)=1$ and $\rho(A,B)=0$ separately. 

\subsection{\boldmath Obstruction to lifting when $\rho\neq0$}

For notational convenience we write $(A_1,A_2)$ instead of $(A,B)$. 
Assuming $\rho(A_1,A_2)=1$, we will show there does not exist a commuting pair $\phi_1,\phi_2\in\Diff(M)$ that induce $A_1,A_2$ on $\pi_1(M)$. 

\subsection*{Lie group reduction}
Since $\Diff(T^d\#\Sigma)\to\SL_d(\Z)$ factors through the mapping class group $\Mod(T^d\#\Sigma)\onto\SL_d(\Z)$, it suffices to show that there is no commuting pair $[\phi_1],[\phi_2]$ in $\Mod(T^d\#\Sigma)$ that lifts $A_1,A_2$. This problem reduces to the problem of lifting $A_1,A_2$ to a commuting pair in the universal cover
\begin{equation}\label{eqn:SL-SES}1\to\Z/2\Z\to\widetilde{\SL_d(\R)}\to\SL_d(\R)\to1\end{equation}
as we now explain.
Let
\[1\to\Z/2\Z\to\widetilde{\SL}_d(\Z)\to\SL_d(\Z)\to1\] be the short exact sequence obtained by pullback of (\ref{eqn:SL-SES}) along the inclusion $\SL_d(\Z)\hookrightarrow\SL_d(\R)$. By \cite[Thm.\ D]{BKKT}, when $\delta(\Sigma)=1$, there is an isomorphism $\Mod(T^d\#\Sigma)\cong K\rtimes\wtil{\SL}_d(\Z)$ (where $K$ is a group whose precise form is not important here), and there is a commutative diagram
\[\begin{xy}
(0,0)*+{K\rtimes\wtil{\SL}_d(\Z)\cong\Mod(T^d\#\Sigma)}="A";
(-20,-15)*+{\wtil{\SL}_d(\Z)}="C";
(20,-15)*+{\SL_d(\Z)}="D";
{\ar "A";"D"}?*!/_3mm/{};
{\ar "C";"D"}?*!/_3mm/{};
{\ar "A";"C"}?*!/_3mm/{};
\end{xy}\]
This implies that if the commuting pair $A_1,A_2$ has a commuting pair of lifts to $\Mod(M)$, then there is also a commuting pair of lifts to $\wtil{\SL}_d(\Z)<\widetilde{\SL_d(\R)}$.

Thus to prove Theorem \ref{thm:nielsen}, it suffices to show that for every pair of lifts $\wtil A_1,\wtil A_2\in\wtil{\SL_d(\R)}$, the commutator $[\wtil A_1,\wtil A_2]$ is nontrivial. Since the extension (\ref{eqn:SL-SES}) is central, it suffices to show this for a single pair of lifts.

To proceed, we use the definition of the universal cover as a set of paths, so choosing lifts amounts to choosing paths from $A_i$ to the identity in $\SL_d(\R)$. With this viewpoint, the commutator $[\wtil A_1,\wtil A_2]$ defines loop in $\SL_d(\R)$ based at the identity, whose homotopy class is an element of the kernel group $\Z/2\Z\cong\pi_1\big(\SL_d(\R)\big)$, and we want to show this homotopy class is nontrivial.

By Theorem \ref{thm:path-components} and its proof, if $\rho(A_1,A_2)=1$, then $(A_1,A_2)$ is connected by a path of commuting pairs in $\SL_d(\R)$ to $(D_1,D_2)$ as defined in \ref{eqn:diagonal-pair}. Therefore, it suffices to compute the obstruction for $(D_1,D_2)$.

\subsection*{Computing the obstruction}

Consider the paths $\eta_i(t)$, $0\le t\le1$, in $\SO(3)$ defined by
\beq
\eta_1(t)=\begin{pmatrix}
\cos(\pi t) & -\sin(\pi t) &0\\
\sin(\pi t) & \cos(\pi t) & 0\\
0 & 0 & 1
\end{pmatrix}
\text{\ \ \ and\ \ \ }
\eta_2(t)=\begin{pmatrix}
1 & 0 & 0\\
0&\cos(\pi t) & -\sin(\pi t)\\
0&\sin(\pi t) & \cos(\pi t) \\
\end{pmatrix}.
\eeq
These define paths from $D_1$ and $D_2$, respectively, to the identity. Note that these are \emph{not} paths of commuting matrices.

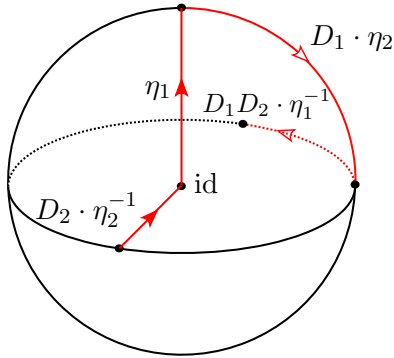
\begin{figure}[h!!]
\labellist
\pinlabel $\text{id}$ at 75 65
\pinlabel $\eta_1$ at 57 100
\pinlabel $D_1\cdot\eta_2$ at 130 120
\pinlabel \small $D_1D_2\cdot\eta_1^{-1}$ at 98 95
\pinlabel $D_2\cdot\eta_2^{-1}$ at 30 55
\endlabellist
\centering
\includegraphics[scale=1]{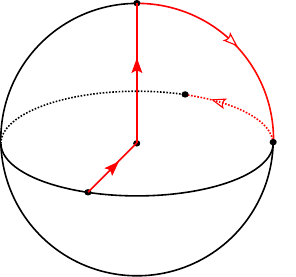}
\caption{Loop homotopic to $[\eta_1,\eta_2]$ in $\SO(3)\cong\R P^3$, viewed as the quotient of the unit 3-ball by the antipodal map on its boundary. A point $v$ in the ball corresponds to the rotation with axis $v$ and angle $|v|\pi$ (counterclockwise according to the right-hand rule). The pictured loop is homotopically nontrivial.}
\label{fig:loop}
\end{figure}

Now we compute the homotopy class of the commutator $[\eta_1,\eta_2]$. Recall that the multiplication in $\wtil{\SL_d(\R)}$ of two paths $\lambda(t),\mu(t)$ in $\SL_d(\R)$ based at the identity is the pointwise product path $t\mapsto \lambda(t)\cdot\mu(t)$ (this holds in any Lie group).
It is helpful to recall that the pointwise product of paths $\lambda,\mu$ is homotopic to the concatenation $\lambda*(\lambda(1)\cdot\eta)$ of $\lambda$ with the path $t\mapsto\lambda(1)\cdot\eta(t)$ (again this holds in any Lie group). Consequently, the path $t\mapsto \eta_1(t)\cdot \eta_2(t)\cdot \eta_1(t)^{-1}\cdot\eta_2(t)^{-1}$ is homotopic to the concatenation of paths
\beq
\eta_1*(D_1\cdot\eta_2)*(D_1D_2\cdot \eta_1^{-1})*(D_1D_2D_1^{-1}\cdot\eta_2^{-1}).
\eeq
Note that $D_1D_2D_1^{-1}=D_2$. One can compute directly that this loop represents a generator of $\pi_1\big(\SO(3)\big)\cong\Z/2\Z$. A picture of this path is given in Figure \ref{fig:loop}.

\subsection{\boldmath Constructing lifts when $\rho=0$}
Assuming $\rho(A_1,A_2)=0$, we construct a splitting of $\ell:\Diff^+(T^d\#\Sigma)\to \SL_d(\Z)$ over $G=\langle A_1,A_2\rangle\cong\Z^2$, i.e.\ we construct an action of $G$ on $T^d\#\Sigma$ that induces the desired splitting.

The standard action of $G<\SL_d(\Z)$ on $T^d=\mathbb R^d/\mathbb Z^d$ fixes the origin $o\in T^d$. By blowing up, we obtain a smooth action of $G$ on a manifold $T'$ diffeomorphic to $T^d\setminus \text{int}(D^d)$ such that the action on the boundary $\partial T'\cong S^{d-1}$ agrees with the action of $G$ on the projectivized tangent space $\mathbb P^+(T_oT^d)\cong S^{d-1}$. In particular, up to this identification, the action of $G$ on $\partial T'$ factors through the standard action of $\SL_d(\mathbb R)$ on $S^{d-1}$ (viewed as rays through the origin).

Now we use the paths constructed in the proof of Theorem \ref{thm:path-components}. In the case $\rho(A_1,A_2)=0$, the concatenation of paths $\# 1$  through $\# 5$ is a path of commuting matrices starting at $(A_1,A_2)$ and ending at the identity matrix. Viewing the path as an isotopy of commuting diffeomorphisms of $S^{d-1}$, we attach a collar to $\partial T'$ on which we perform this isotopy (we can arrange that the path is smooth and locally constant near the boundary, so that the resulting action is smooth). Now we have an action of $G$ on a manifold $T''$ that is diffeomorphic to $T'$, such that $G$ acts trivially near $\partial T''$. We then glue $T''$ and $D^d$ along their boundary using a diffeomorphism representing $\Sigma\in\Theta_d\cong\pi_0\Diff^+(S^{d-1})$. Extending the action on $T''$ to $T^d\#\Sigma\cong T''\cup D^d$ by the identity on $D^d$, defines the desired action of $G$ on $T^d\#\Sigma$. 
This concludes the proof of Theorem \ref{thm:nielsen}.\qed
 
\section{Homotopy Anosov actions (Proof of Theorem \ref{thm:anosov})}\label{sec:anosov}

To prove Theorem \ref{thm:anosov}, fix $d$ such that there exists $\Sigma\in\Theta_d$ satisfying $\delta(\Sigma)=1$. 

First we construct (infinitely many conjugacy classes of pairs of) commuting hyperbolic matrices $A_1,A_2\in\SL_d(\Z)$ such that $\rho(A_1,A_2)=1$. 
(Recall that a matrix in $\SL_d(\Z)$ is called hyperbolic if it has no eigenvalues on the unit circle.) 
Write $d=n+3$ and define $A_i$ to be a block diagonal matrix $\left(\begin{smallmatrix}B_i&\\&C_i\end{smallmatrix}\right)$, where $B_1,B_2\in\SL_3(\Z)$ and $C_1,C_2\in\SL_n(\Z)$ are commuting pairs of hyperbolic matrices as defined in the following paragraphs.

Choose $B_1,B_2\in\SL_3(\Z)$ that are conjugate in $\SL_3(\R)$ to diagonal matrices of the form
\begin{equation}\label{eqn:matrices}
\begin{pmatrix}
\lambda_1 & 0 &0\\
0 & \lambda_2 & 0\\
0 & 0 & \frac{1}{\lambda_1\lambda_2}
\end{pmatrix}
\text{\ \ \ and\ \ \ }
\begin{pmatrix}
\frac{1}{\mu_1\mu_2} & 0 & 0\\
0 & \mu_1 & 0\\
0 & 0 & \mu_2
\end{pmatrix}
\end{equation}
respectively, where $\lambda_1,\lambda_2,\mu_1,\mu_2$ are all negative and different from $-1$. As an explicit example, consider the polynomial $\xi=x^3 + x^2 - 2 x - 1$. The totally real cubic field $K=\Q[x]/(\xi)$ has discriminant 49 (the smallest possible among totally real cubic fields). Fixing a root $\alpha$ of $\xi$ in $K$, the group of units $\mathcal O_K^\times$, modulo its torsion subgroup (which is isomorphic to $\Z/2\Z$, generated by $-1$), is freely generated by $\epsilon_1:=\alpha^2+\alpha-1$ and $\epsilon_2:=-\alpha^2+2$. The action of the units $-\epsilon_1$ and $\epsilon_1\epsilon_2$ on the ring of integers $\mathcal O_K$ with the basis $\mathcal O_K\cong \Z\{1,\alpha,\alpha^2\}$ gives matrices as in (\ref{eqn:matrices}). For the claims about this number field, compare with \cite[\S B.4]{cohen}.

To choose $C_1,C_2\in\SL_n(\Z)$, we recall that for each $n\ge3$, there exists a subgroup $\Z^2<\SL_n(\Z)$ generated by hyperbolic matrices $C_1,C_2$ such that all eigenvalues of $C_1$ and $C_2$ are real and positive. Indeed, let $K/\Q$ be a degree $n$ totally real number field. Choose linearly independent units $\alpha_1,\alpha_2\in\mathcal O_K^\times$, and let $C_i$ be the matrix for multiplication by $\alpha_i$ on $\mathcal O_K\cong\Z^n$ (with respect to any basis). Since the Galois conjugates of the $\alpha_i$ are real and not equal to $\pm1$, they do not lie on the unit circle, so the matrices $C_i$ are hyperbolic. Furthermore, after replacing $\alpha_i$ by $\alpha_i^2$, we can ensure that the eigenvalues of $C_i$ are positive.

Now define $G=\langle A_1,A_2\rangle$ and consider the Anosov action of $G$ on $T^d$. This action is not the linearization of any smooth action on  $ T^d\#\Sigma$ because if it were, then $\Diff^+(T^d\#\Sigma)\to\SL_d(\mathbb Z)$ would split over $G$, contradicting the Theorem \ref{thm:nielsen}. Here we use the fact that $\delta(\Sigma)=1$ by assumption and $\rho(A_1,A_2)=1$ by construction. This completes the proof. 
\qed

\section{Bordism of group actions (Theorem \ref{thm:bordism})}\label{sec:bordism}

Fix a commuting pair $A,B$ in $\SL_d(\R)$ and fix $f\in\Diff^+(S^{d-1})$. We view $A,B$ as diffeomorphisms of $S^{d-1}$, via the action on rays in $\R^d$; since $\SL_d(\R)$ is path-connected, these diffeomorphisms are isotopic to the identity. Consider the subgroup
$G=G(f,A,B)$ of $\Diff_0(S^{d-1})$ generated by $fAf^{-1}$ and $fBf^{-1}$. Here $\Diff_0$ denotes diffeomorphisms isotopic to the identity.

The following statement is a more precise form of Theorem \ref{thm:bordism}. We identify $\pi_0\Diff^+(S^{d-1})$ with $\Theta_d$ (via the clutching construction). In particular, for each isotopy class $[f]\in\pi_0\Diff^+(S^{d-1})$, there is an associated element $\eta\cdot[f]\in\Theta_{d+1}$ given by the Milnor--Munkres--Novikov pairing. 
\begin{thm}\label{thm:no-extend}
Let $A,B\in\SL_d(\R)$ a pair of commuting matrices, $f\in\Diff^+(S^{d-1})$, and $G=G(A,B,f)<\Diff_0(S^{d-1})$ as above. If $\rho(A,B)=1$ and $\eta\cdot[f]\in\Theta_{d+1}$ is nonzero, then there is no action of $G$ on $D^d$ that extends the given action of $G$ on $S^{d-1}$.
\end{thm}

In terms of examples, we could take $(A,B)=(D_1,D_2)$ from (\ref{eqn:diagonal-pair}) to get examples with $G\cong\Z/2\Z\times \Z/2\Z$. We can also easily find $A,B$ so that $G\cong\Z^2$; for example

\[A=
\left(\begin{smallmatrix}
-1&\\
&-1\\
&&1\\
&&&1&1\\
&&&&1\\
&&&&&1\\
&&&&&&1
\end{smallmatrix}\right)\>\>\text{ and }B=\left(\begin{smallmatrix}
-1&\\
&1\\
&&-1\\
&&&1\\
&&&&1\\
&&&&&1&1\\
&&&&&&1
\end{smallmatrix}\right).\]
We also recall that the function $\eta: \Theta_d\to\Theta_{d+1}$ is nontrivial for infinitely many values of $d$; see \cite{BKKT}.

Before proving Theorem \ref{thm:no-extend}, we give two explicit descriptions of the Milnor--Munkres--Novikov pairing $\eta\cdot[f]$ for $[f]\in\pi_0\Diff^+(S^{d-1})$ that will be used in the proof. For reference, see \cite[\S2]{munkres} and \cite[\S7]{levine}.

\subsubsection*{MMN pairing: first version}\label{sec:MMN1}
Let $d\geq 3$. Write $r_t$ for a loop in $\SO(d)$ based at the identity that generates $\pi_1(\SO(d))\cong\Z/2\Z$. Given a diffeomorphism $f\in\Diff_{\partial}(D^{d-1})$ which is the identity near the boundary, define a diffeomorphism $g$ of $D^d=[0,1]\times D^{d-1}$ rel boundary by
\[g_f(t,x)= \big(t,[r_t,f](x)\big).\]
Note that we can write $g_f$ as the commutator of diffeomorphisms $R,F$ where $R(t,x)=(t,r_t(x))$ and $F(t,x)=(t,f(x))$.

\subsubsection*{MMN pairing: second version}\label{sec:MMN2}
As a variation of the previous construction, take $r_t$ as above, fix $f\in\Diff^+(S^{d-1})$, and consider the diffeomorphism of $[0,1]\times S^{d-1}$ given by $FR^{-1}F^{-1}$, i.e.\ the map $(t,x)\mapsto \big(t, (f\circ r_t^{-1}\circ f^{-1})(x)\big)$. This diffeomorphism is the identity on the boundary, so fixing a (standard) embedding of $[0,1]\times S^{d-1}$ in $D^d\subset S^d$, we can extend by the identity to get a diffeomorphism that represents the same isotopy class as $g_f$. To see that the isotopy class is the same, it helps to observe that the diffeomorphism $R$ of $[0,1]\times S^{d-1}$ extends to $[0,1]\times D^d$, and using this, we can conclude that once we extend $R$ by the identity to $S^d$, we get a diffeomorphism that extends to $D^{d+1}$ and hence is (pseudo-)isotopically trivial. Then $FR^{-1}F^{-1}$ and $RFR^{-1}F^{-1}$ are isotopic.

\begin{proof}[Proof of Theorem \ref{thm:no-extend}]
Extending the given action of $G$ on $S^{d-1}$ to a smooth action on $D^d$ amounts to solving the following lifting problem.
\beq
\xymatrix{
1\ar[r]&\Diff_{\partial}(D^d)\ar[r]&\Diff_0(D^d)\ar[r]& \Diff_0(S^{d-1})\ar[r]&1\\
&&&G\ar@{^{(}->}[u]\ar@{.>}[lu]
}
\eeq
Pulling back this sequence along $G\hookrightarrow\Diff_0(S^{d-1})$ and taking path components, we obtain an extension
\begin{equation}\label{eqn:extend-disk-extension}
1\to \Theta_{d+1} \to E\to G \to 1,
\end{equation}
where we use the isomorphisms $\Theta_{d+1}\cong \pi_0\Diff^+(S^{d}\big)\cong \pi_0 \Diff_{\partial}(D^d)$. 

Observe that the obstruction to splitting the extension \eqref{eqn:extend-disk-extension} gives an obstruction to the original lifting problem. The extension \eqref{eqn:extend-disk-extension} is central because $\Theta_{d+1}$ is abelian and because any $g\in G<\Diff^+(S^{d-1})$ has a lift to $\Diff_0(D^d)$ that is supported near $\partial D^d$, while every element of $\Diff_\partial (D^d)$ can be isotoped to have support away from $\partial D^d$. 
Because the extension \eqref{eqn:extend-disk-extension} is central, the obstruction to splitting is the commutator of the lifts of the generators of $G$ \cite[\S IV.3.\ Exercise 8]{Brown}. The diffeomorphisms $fAf^{-1}$ and $fBf^{-1}$ of $S^{d-1}$ can be extended radially using isotopies $A_t$ and $B_t$ from $A$ and $B$ to the identity in $\SL_d(\R)$. Using these isotopies, we extend $fAf^{-1}$ and $fBf^{-1}$ radially on a collar neighborhood of $\partial D^d$. This defines elements $a,b\in \Diff_0(\partial D^d)$ and we write $\alpha,\beta$ for the corresponding elements of $E$; these are lifts of $fAf^{-1},fBf^{-1}\in G$, respectively. Our goal is to show $[\alpha,\beta]\neq0$, or equivalently that $[a,b]\in\Diff_\partial(D^d)$ is not isotopic to the identity. 

The commutator $[a,b]$ is supported on a collar $[0,1]\times S^{d-1}$, where it is given by the formula
\[(t,x)\mapsto \left(t, \big(f\left[A_t,B_t\right]f^{-1}\big)(x)\right).\]
To show that $[a,b]$ is not isotopic to the identity, 
recall from the proof of Theorem \ref{thm:nielsen} that if $\rho(A,B)=1$, then the path $r_t:=[A_t,B_t]$ defines a nontrivial loop in $\pi_1\big(\SL_d(\R)\big)$. Now by the preceding discussion of the MMN pairing, the diffeomorphism $(t,x)\mapsto \big(t, (f\circ r_t^{-1}\circ f^{-1})(x)\big)$ represents the isotopy class $\eta\cdot[f]\in\Theta_{d+1}$. Since by assumption $\eta\cdot[f]\neq0$, it follows that the isotopy class of $[a,b]$ is nontrivial, as desired. Consequently, the extension (\ref{eqn:extend-disk-extension}) does not split.
\end{proof}

\bibliographystyle{amsalpha}
\bibliography{refs}

Mauricio Bustamante\\
Departamento de Matem\'aticas, Pontificia Universidad Cat\'olica de Chile\\
\texttt{mauricio.bustamante@uc.cl}

Bena Tshishiku\\
Department of Mathematics, Brown University\\
\texttt{bena\_tshishiku@brown.edu}

\end{document}